\documentclass{mcom-l}
\usepackage[top=3.54cm,bottom=3.54cm,left=2.94cm,right=2.94cm ]{geometry}

\usepackage{amssymb}

\usepackage{graphicx}

\usepackage{color}
\usepackage[cmtip,all]{xy}

\usepackage{graphicx}
\usepackage{mathtools}
\usepackage{pifont}
\usepackage{tikz}
\usetikzlibrary{arrows.meta}
\usetikzlibrary{positioning}
\allowdisplaybreaks

\newtheorem{theorem}{Theorem}[section]

\theoremstyle{definition}

\theoremstyle{remark}

\numberwithin{equation}{section}

\renewcommand{\theequation}{\arabic{section}.\arabic{equation}}
\renewcommand{\thetable}{\arabic{section}.\arabic{table}}

\newcommand{\be}{\begin{equation}}
\newcommand{\ee}{\end{equation}}
\newcommand{\ba}{\begin{array}}
\newcommand{\ea}{\end{array}}
\newcommand{\bea}{\begin{eqnarray}}
\newcommand{\eea}{\end{eqnarray}}
\newcommand{\beas}{\begin{eqnarray*}}
\newcommand{\eeas}{\end{eqnarray*}}

\begin{document}

\title[]{Characterizations of the generalized inverse Gaussian, asymmetric Laplace, and shifted (truncated) exponential laws via independence properties}


\author{Kevin B. Bao}
\address{Cornell University, Ithaca, NY, USA}
\email{kbb37@cornell.edu}

\author{Christian Noack}
\address{Department of Mathematics, Cornell University, Ithaca, NY, USA}
\email{noack@cornell.edu}

\subjclass[2020]{}

\date{}

\dedicatory{}

\begin{abstract}
	We prove three new characterizations of the generalized inverse Gaussian (GIG), asymmetric Laplace (AL), shifted exponential (sExp) and shifted truncated exponential (stExp) distributions in terms of non-trivial independence preserving transformations, which were conjectured by Croydon and Sasada in \cite{CS1}. We do this under the assumptions of absolute continuity and mild regularity conditions on the densities.  

    Croydon and Sasada \cite{CS2} use these independence preserving transformations to analyze statistical mechanical models which display KPZ behavior. Our characterizations show the integrability of these models only holds for these four specific distributions in the absolutely continuous setting.
\end{abstract}

\maketitle

\section{Introduction}\setcounter{equation}{0}

\subsection{Background}
Several important distributions preserve independence under certain transformations. Perhaps the most classic example is the normal distribution. Kac \cite{K} and Bernstein \cite{B} both showed that given two non-degenerate real independent random variables $X$ and $Y$, the two random variables $(U,V) \coloneqq \left(\frac{X + Y}{2}, \frac{X - Y}{2} \right)$ are independent if and only if $X$ and $Y$ are normally distributed with the same variance.  This in turn implies that $U$ and $V$ are both normal with the same variance as well.  Note that this transformation is an involution, a property which all transformations in this paper will possess.

Similar independence properties have also been proven for the gamma distribution in Lukacs \cite{L}, the geometric and exponential distributions in Crawford \cite{C},  the beta distribution in Seshadri and Wesołowski \cite{SW_2003}, and the product of generalized inverse Gaussian and gamma distributions in Matsumoto and Yor \cite{MY} and Seshadri and Wesołowski \cite{SW_2001}.

The characterizations of Lukacs \cite{L} and Seshadri and Wesołowski \cite{SW_2003}  were used in Chaumont and Noack \cite{CN} to characterize  stationarity in $1+1$ dimensional lattice directed polymer models, which in turn characterize the Burke property in the same setting. While many statistical mechanical models are conjectured to belong to the KPZ universality class, computations quickly become intractable without the presence of some sort of integrability structure such as stationarity. In fact, at the moment, only exactly solvable models have allowed for KPZ-type computations, highlighting how useful independence preserving transformations can be. The characterization of independence preserving distributions can therefore aid in the search for new exactly solvable models displaying KPZ behavior.

In this paper, we prove characterizations via independence preserving transformations for the generalized inverse Gaussian (GIG), the asymmetric Laplace (AL), the shifted exponential (sExp), and the shifted truncated exponential (stExp) distributions, all of which were conjectured in Croydon and Sasada \cite{CS1}.

\subsection{Main Results}
We now recall the definitions of the distributions followed by the theorems that characterize them.

{\bf Generalized inverse Gaussian (GIG) distribution:} For $\lambda \in \mathbb{R}$, $c_1,c_2 \in (0, \infty)$, the generalized inverse Gaussian distribution with parameters $(\lambda, c_1, c_2)$, which we denote GIG$(\lambda,c_1,c_2)$, has density
\[
    \dfrac{1}{Z} x^{\lambda - 1} e^{-c_1x-c_2x^{-1}} \mathbf{1}_{(0,\infty)}(x), \qquad x \in \mathbb{R},
\]
where $Z = \frac{2K_\lambda(2\sqrt{c_1c_2})}{\left(\frac{c_1}{c_2}\right)^{\lambda/2}}$ is the normalizing constant and $K_\lambda$ is the modified Bessel function of the second kind with parameter $\lambda$.

\begin{theorem}\label{thm:gig}
    Let $\alpha,\beta \geq 0$ with $\alpha \neq \beta$ and $F_1:(0,\infty)^2 \rightarrow (0,\infty)^2$ be the involution given by
    \be\label{eq:gigtrans}
        F_1(x,y) = \left(\dfrac{y(1 + \beta xy)}{1 + \alpha xy}, \dfrac{x(1 + \alpha xy)}{1 + \beta xy}\right).
    \ee
    Let X and Y be $(0,\infty)$-valued independent random variables with twice-differentiable densities that are strictly positive throughout $(0,\infty)$. Then, the random variables $(U,V)\coloneqq$ $F_1(X,Y)$ are independent if and only if there exist $\lambda\in \mathbb{R}$ and $c_1,c_2 > 0$ such that
    \[
        X \sim \textup{GIG}(\lambda,c_1\alpha,c_2) \quad and \quad Y \sim \textup{GIG}(\lambda,c_2\beta,c_1),
    \]
    in which case, U $\sim$ \textup{GIG}$(\lambda,c_2\alpha,c_1)$ and V $\sim$ \textup{GIG}$(\lambda,c_1\beta,c_2)$. Hence, if moreover (U,V) has the same distribution as (X,Y), then X $\sim$ \textup{GIG}$(\lambda,c\alpha,c)$ and Y $\sim$ \textup{GIG}$(\lambda,c\beta,c)$ for some $\lambda\in \mathbb{R}$ and $c > 0$.
\end{theorem}

\textbf{Remark:} A special case of this theorem with $\alpha=1$, $\beta=0$ has already been proven by Letac and Weso\l owski \cite[Theorem 4.1]{LW} where they used no density assumptions.

Now recall the definition of the asymmetric Laplace (AL) distribution.

{\bf Asymmetric Laplace (AL) distribution:} For $\lambda_1,\lambda_2 \in (0,\infty)$, the asymmetric Laplace distribution with parameters $(\lambda_1,\lambda_2)$, which we denote AL$(\lambda_1,\lambda_2)$, has density
\[
    \dfrac{1}{Z} \left( e^{-\lambda_1 x} \mathbf{1}_{[0,\infty)}(x) + e^{\lambda_2 x} \mathbf{1}_{(-\infty,0)}(x) \right), \qquad x \in \mathbb{R},
\]
where $Z = \frac{1}{\lambda_1} + \frac{1}{\lambda_2} \left(\text{or}\ \frac{1}{Z} = \frac{\lambda_1\lambda_2}{\lambda_1 + \lambda_2} \right)$.

\begin{theorem}\label{thm:al}
	Let $F_2:\mathbb{R}^2 \rightarrow \mathbb{R}^2$ be the involution given by
    \be\label{eq:altrans}
        F_2(x,y) = (\min\{x,0\} - y, \min\{x,y,0\} - x - y).
    \ee
    Let X and Y be $\mathbb{R}$-valued independent random variables with densities that never vanish and are twice-differentiable on $(-\infty,0)\cup(0,\infty)$. It is then the case that (U,V) $\coloneqq F_2(X,Y)$ are independent if and only if there exist $p,q,r > 0$ such that
    \[
        X \sim \textup{AL}(p,q) \quad and \quad Y \sim \textup{AL}(p+q,r),
    \]
    in which case, U $\sim$ \textup{AL}(r,q) and V $\sim$ \textup{AL}(q+r,p). Hence, if moreover (U,V) has the same distribution as (X,Y), then X $\sim$ \textup{AL}(p,q) and Y $\sim$ \textup{AL}(p+q,p).
\end{theorem}

Notice that Theorem 1.2 has striking similarities with Wesołowski \cite{W}, who showed that given two non-degenerate $(0,1)$-valued independent random variables $X,Y$, the random variables $(U,V) \coloneqq G(X,Y)$ under the transformation
\[
    G(x,y) = \left(\dfrac{1-y}{1-xy},1-xy \right)
\]
are independent if and only if $(X,Y) \sim (\beta_{p,q},\beta_{p+q,r})$ where $\beta_{p,q}$ denotes the beta distribution with shape parameters $p,q$. Theorem 1.2 has a similar condition $(X,Y) \sim (\textup{AL}_{p,q},\textup{AL}_{p+q,r})$ where $\textup{AL}_{p,q}$ denotes the asymmetric Laplace (AL) distribution with parameters $p,q$.

Finally, we recall the definitions of the shifted exponential (sExp) and shifted truncated exponential (stExp) distributions. The definitions are then followed by Theorem 1.3 which characterizes them.

{\bf Shifted exponential distribution:} For $\lambda > 0,c \in \mathbb{R}$, the shifted exponential distribution with parameters $(\lambda,c)$, which we denote sExp$(\lambda,c)$, has density
\[
    \dfrac{1}{Z} e^{-\lambda x} \mathbf{1}_{[c,\infty)}(x), \qquad x \in \mathbb{R},
\]
where $Z = \frac{1}{\lambda} e^{-\lambda c}$.

{\bf Shifted truncated exponential distribution:} For $\lambda > 0,c_1,c_2 \in \mathbb{R}$ with $c_1 < c_2$, the shifted truncated exponential distribution with parameters $(\lambda,c_1,c_2)$, which we denote stExp$(\lambda,c_1,c_2)$, has density
\[
    \dfrac{1}{Z} e^{-\lambda x} \mathbf{1}_{[c_1,c_2]}(x), \qquad x \in \mathbb{R},
\]
where $Z = \frac{1}{\lambda} \left(e^{-\lambda c_1} - e^{-\lambda c_2}\right)$.

\begin{theorem}\label{thm:se}
Let $c_1, c_2>0$ and define
\begin{equation}\label{eq:sExptrans}
    F_3(x,y) = (\min\{-x,y\}, y + x - \min\{-x,y\}).
\end{equation}
Then $F_3: \mathbb{R}^2 \rightarrow \mathbb{R}^2$ is an involution, while $F_3: [-c_1, c_2] \times [-c_2, \infty) \rightarrow [-c_2, c_1] \times [-c_1, \infty)$ is a bijection. Let X and Y be $\mathbb{R}$-valued independent random variables with densities  satisfying $f_X(x) > 0$ when $x \in [-c_1, c_2]$ with $f_X(x) = 0$ otherwise, and $f_Y(y) > 0$ when $y \in [-c_2, \infty)$ with $f_Y(y) = 0$ otherwise. Moreover, assume the densities are twice-differentiable inside their respective supports. It is then the case that (U,V) $\coloneqq F_3(X,Y)$ are independent if and only if there exists $\lambda > 0$ such that
    \[
        X \sim \textup{stExp}(\lambda, -c_1, c_2) \quad and \quad Y \sim \textup{sExp}(\lambda, -c_2),
    \]
    in which case, U $\sim$ \textup{stExp}$(\lambda, -c_2, c_1)$ and V $\sim$ \textup{sExp}$(\lambda, -c_1)$. Hence, if moreover (U,V) has the same distribution as (X,Y), then X $\sim$ \textup{stExp}$(\lambda,-c,c)$ and Y $\sim$ \textup{sExp}$(\lambda,-c)$ for some $c > 0$. Note that $F_3: [-c, c] \times [-c,\infty) \rightarrow [-c, c] \times [-c,\infty)$ is an involution.
\end{theorem}

\subsection{Structure of the paper}
Each of our main results correspond to a conjecture in  Croydon and Sasada \cite{CS1}. Theorem 1.1 solves Conjecture 8.6, Theorem 1.2 solves Conjecture 8.15, and Theorem 1.3 solves Conjecture 8.10.

We first prove Theorem 1.2 and 1.3 and then prove Theorem 1.1. We save the proof of Theorem 1.1 for the end because it is the most involved of the three. In Section 2, we give the proof for Theorem 1.2. We then give the proof for Theorem 1.3 in Section 3. Finally, we prove Theorem 1.1 in Section 4. All three theorems are proven using methods from Chaumont and Noack in \cite{CN}.

\subsection{Acknowledgements}
The authors would like to thank Timo Sepp{\"a}l{\"a}inen and Philippe Sosoe for their valuable insights.

\section{Proof of Theorem 1.2}
We begin by partitioning $\mathbb{R}^2$ into sections so that the minimum functions in \eqref{eq:altrans} may be simplified. Define connected open subsets of $\mathbb{R}^2$ by
\begin{align*}
    &\text{\ding{192}} :=\{{(x,y):x>0, y>0}\}, &\text{\ding{193}} :=\{{(x,y):x<0, y>0}\}, \qquad &\text{\ding{194}} :=\{{(x,y):x < y < 0}\},\\ &\text{\ding{195}} :=\{{(x,y):y < x < 0}\}, &\text{\ding{196}} :=\{{(x,y):x>0, y<0}\}. \qquad &
\end{align*}
See Figure 1 below.

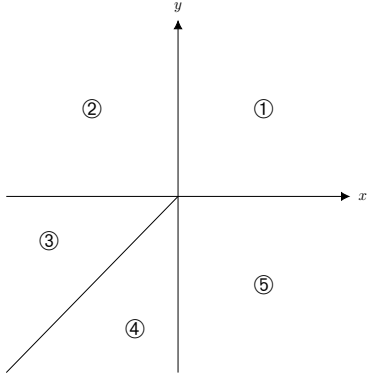
\begin{figure}[ht]
\centering
\resizebox{5cm}{5cm}{
    \begin{tikzpicture}
        \draw [{Latex[length = 2mm, width = 2mm]}-](0,4) -- (0,-4);
        \draw [{Latex[length = 2mm, width = 2mm]}-](4,0) -- (-4,0);
        \draw (-4,-4) -- (0,0);
        \node at (2,2) {\huge{\ding{192}}};
        \node at (-2,2) {\huge{\ding{193}}};
        \node at (-3,-1) {\huge{\ding{194}}};
        \node at (-1,-3) {\huge{\ding{195}}};
        \node at (2,-2) {\huge{\ding{196}}};
        \node at (4.3,0) {$x$};
        \node at (0,4.3) {$y$};
    \end{tikzpicture}
}
\caption{Partition of $\mathbb{R}^2$}
\label{fig:1}
\end{figure}

Recall that the transformation in this theorem is
\[
    F_2(x,y) = (\min\{x,0\} - y, \min\{x,y,0\} - x - y)=:(u,v),
\]
which maps \ding{192} $\leftrightarrow$ \ding{195}, \ding{193} $\leftrightarrow$ \ding{194}, and \ding{196} $\leftrightarrow$ \ding{196}. Therefore, defining the set $O := $ \ding{192} $\cup$ \ding{193} $\cup$ \ding{194} $\cup$ \ding{195} $\cup$ \ding{196}, we see that $F_2$ is a smooth involution on the open set $O$. Notice that $\mathbb{R}^2 \setminus O$ has Lebesgue measure $0$ and therefore has no effect on the distributions of the pairs of random variables $(X,Y)$ and $(U,V)$. We will therefore only focus on the joint densities with inputs $(x,y)\in O$, which in turn implies $(u,v)=F_2(x,y)\in O$.

Our approach begins with the joint densities of $(U,V)$ and $(X,Y)$, which are related in the following manner
\begin{equation}
    f_{UV}(u,v) = |J| f_{XY}(x,y)   \qquad \text{ for } (x,y) \in O\label{eq: invol density 1}
\end{equation}
where $f_A(a)$ denotes the probability distribution function of  $A$, and $J$ is the Jacobian determinant of the involution $F_2$.

One can  check (as we do later) that $J\equiv -1$ throughout \ding{192}, \ding{193}, \ding{194}, \ding{195}, \ding{196} and therefore throughout O.

\begin{proof}
The \textit{if} part follows from an easy computation using \eqref{eq: invol density 1}.
     We now prove the \textit{only if} part.  Since $U,V$ are mutually independent and $X,Y$ are mutually independent, \eqref{eq: invol density 1} simplifies to
    \begin{equation}
        f_U(u) f_V(v) = f_X(x) f_Y(y) \qquad \text{ for } (x,y)\in O.\label{eq:al1}
    \end{equation}
Since neither $f_X$ nor $f_Y$ vanish on $(-\infty,0)\cup (0,\infty)$, \eqref{eq:al1} implies that neither $f_U$ nor $f_V$ vanish on $(-\infty, 0) \cup (0, \infty)$. Since $F_2$ is a smooth involution on $O$, and $f_X$ and $f_Y$ are both twice-differentiable on $(-\infty,0)\cup (0,\infty)$, the same follows for $f_U$ and $f_V$.

    In section \ding{192}, where $x > 0$ and $y > 0$, the transformation simplifies to
    \[
        (u,v) = F_2(x,y) = (-y,-x-y)
    \]
    which has Jacobian determinant $J = \begin{vmatrix} 0 & -1 \\ -1 & 1 \end{vmatrix} = -1$.
    Plugging the transformation into \eqref{eq:al1}, we get
    \begin{equation}\label{eq:al2}
        f_U(-y) f_V(-x-y) = f_X(x) f_Y(y) \qquad \text{ for all } x,y>0.
    \end{equation}
    Taking the logarithm of both sides, and setting $r_A(a) = \log f_A(a)$ gives
    \begin{equation}\label{eq:al3}
        r_U(-y) + r_V(-x-y) = r_X(x) + r_Y(y).
    \end{equation}
    Differentiating \eqref{eq:al3} with respect to $x$ yields
    \begin{equation}\label{eq:al4}
        -{r_V}'(-x-y) = {r_X}'(x).
    \end{equation}
    Now differentiate \eqref{eq:al4} with respect to $y$ to obtain
    \begin{equation}\label{eq:al5}
        0 = {r_V}''(-x-y) = {r_V}''(v).
    \end{equation}
    The second equality holds for all $v<0$ because in section \ding{192}, $v=-x-y < 0$. Thus $r_V$ must be linear on $(-\infty,0)$ and there must exist some real constants $a_1, a_2$ such that
    \begin{equation}\label{eq:al6}
        r_V(v) = a_1 v + a_2 \qquad  \text{ for all } v < 0.
    \end{equation}
     This implies
    \begin{equation}\label{eq:al v<0}
        f_V(v) = e^{a_1 v + a_2} \qquad \text{ for all } v < 0.
    \end{equation}
    Substituting the derivative of \eqref{eq:al6} into \eqref{eq:al4} gives ${r_X}'(x) = -a_1$ for all $x>0$, which in turn implies the existence of some real constant $a_3$ such that

    \begin{equation}\label{eq:al x>0}
        f_X(x) = e^{-a_1 x + a_3}  \qquad \text{ for all } x > 0.
    \end{equation}
    We have thus shown that the p.d.f.'s of $X$ and $V$ have exponential forms in the domains $(0, \infty)$ and $(-\infty, 0)$ respectively.
    
    In section \ding{193}, where $x <0$ and $y > 0$, the transformation simplifies to
    \[
        (u,v) = F_2(x,y) = (x-y,-y)
    \]
    whose Jacobian determinant is $J = \begin{vmatrix} 1 & -1 \\ 0 & -1 \end{vmatrix} = -1$.
    Plugging this transformation into \eqref{eq:al1}, we get
    \begin{equation}\label{eq:al7}
        f_U(x-y) f_V(-y) = f_X(x) f_Y(y) \qquad \text{ for all } x<0 \text{ and } y>0.
    \end{equation}
    Taking the logarithm of both sides, setting $r_A(a) = \log f_A(a)$, and differentiating with respect to $x$ gives
    \begin{equation}\label{eq:al9}
        {r_U}'(x-y) = {r_X}'(x).
    \end{equation}
    Now differentiate \eqref{eq:al9} with respect to $y$ to obtain
    \begin{equation}\label{eq:al10}
        0 = -{r_U}''(x-y) = -{r_U}''(u).
    \end{equation}
    The second equality holds for all $u<0$, because in section \ding{193}, $u=x-y< 0$. Thus $r_U$ must be linear throughout $(-\infty,0)$ and there must exist real constants $a_4,a_5$ such that
    \begin{equation}\label{eq:al11}
        r_U(u) = a_4 u + a_5 \qquad \text{ for all } u < 0.
    \end{equation}
    This implies that
    \begin{equation}\label{eq:al u<0}
        f_U(u) = e^{a_4 u + a_5} \qquad \text{ for all } u < 0.
    \end{equation}
    Differentiating \eqref{eq:al11} and substituting into \eqref{eq:al9} gives $ {r_X}'(x) = a_4$ for all $x<0$ which implies the existence of a real constant $a_6$ such that
    \begin{equation}\label{eq:al x<0}
        f_X(x) = e^{a_4 x + a_6} \qquad \text{ for all } x < 0.
    \end{equation}
    We have thus shown that the p.d.f.'s of $X$ and $U$ both have exponential forms in the domain $(-\infty, 0)$.
    
    In section \ding{195}, where $y < x < 0$, the transformation simplifies to
    \[
        (u,v) = F_2(x,y) = (x-y,-x)
    \]
    which has Jacobian determinant $J = \begin{vmatrix} 0 & -1 \\ -1 & -1 \end{vmatrix} = -1$.
    Plugging this transformation into \eqref{eq:al1}, we get
    \begin{equation}\label{eq:al12}
        f_U(x-y) f_V(-x) = f_X(x) f_Y(y) \qquad \text{ for all } y<x<0.
    \end{equation}
    Taking the logarithm of both sides, setting $r_A(a) = \log f_A(a)$, and differentiating with respect to $y$ gives
    \begin{equation}\label{eq:al14}
        -{r_U}'(x-y) = {r_Y}'(y).
    \end{equation}
    Now differentiate \eqref{eq:al14} with respect to $x$ to obtain
    \begin{equation}\label{eq:al15}
        0 = -{r_U}''(x-y) = {r_U}''(u).
    \end{equation}
    The second equality holds for all $u>0$ because in \ding{195}, $u=x-y > 0$. Thus $r_U$ must be linear throughout $(0,\infty)$, meaning there exist real constants $a_7, a_8$ such that
    \begin{equation}\label{eq:al16}
        r_U(u) = -a_7 u + a_8  \qquad \text{ for all } u > 0.
    \end{equation}
    This implies that
    \begin{equation}\label{eq:al u>0}
        f_U(u) = e^{-a_7 u + a_8} \qquad  \text{ for all } u > 0.
    \end{equation}
    Substituting the derivative of \eqref{eq:al16} into \eqref{eq:al14} gives ${r_Y}'(y) = a_7$ for all $y<0$, which in turn implies the existence of a real constant $a_9$ such that
    \begin{equation}\label{eq:al y<0}
        f_Y(y) = e^{a_7 y + a_9}\qquad \text{ for all } y < 0.
    \end{equation}
    We have thus shown that the p.d.f.'s of $Y$ and $U$ have exponential forms in the domains $(-\infty, 0)$ and $(0, \infty)$ respectively.
    
    So far, we have proven that the p.d.f.'s of $X$ and $U$ both have exponential forms in the domain $(-\infty, 0) \cup (0, \infty)$ and the p.d.f's of $Y$ and $V$ both have exponential forms in the domain $(-\infty, 0)$. Now, if we substitute \eqref{eq:al u<0}, \eqref{eq:al v<0}, \eqref{eq:al x>0} into \eqref{eq:al2} and \eqref{eq:al u>0}, \eqref{eq:al x<0}, \eqref{eq:al y<0} into \eqref{eq:al12}, we will see that $Y$ and $V$ also have exponential forms in the domain $(0, \infty)$.
    
    \begin{equation*}
        \begin{dcases}
            e^{-a_4 y + a_5} e^{-a_1(x+y) + a_2} = e^{-a_1x + a_3} f_Y(y) \\
            e^{-a_7(x-y) + a_8} f_V(-x) = e^{a_4 x + a_6} e^{a_7 y + a_9}
        \end{dcases}
    \end{equation*}
    
    Hence, we have
    \begin{equation}\label{eq:al y>0}
        f_Y(y) = e^{-(a_1 + a_4)y + a_2 + a_5 - a_3}\qquad \text{ for all } y > 0,
    \end{equation}
    \begin{equation}\label{eq:al v>0}
        f_V(v) = e^{-(a_4 + a_7)v + a_6 + a_9 - a_8}\qquad \text{ for all } v > 0.
    \end{equation}

    Now that we have shown all the random variables have exponential forms everywhere, we proceed to show that for each of the densities $f_X, f_Y, f_U, f_V$ the respective constants coefficients in the domains $(-\infty, 0)$ and $(0, \infty)$ are equal. Substituting \eqref{eq:al u<0}, \eqref{eq:al v<0}, \eqref{eq:al x<0}, and \eqref{eq:al y>0} into \eqref{eq:al7} gives
    \begin{equation*}
        e^{a_4 (x-y) + a_5} e^{-a_1 y + a_2} = e^{a_4 x + a_6} e^{-(a_1 + a_4)y + a_2 + a_5 - a_3},
    \end{equation*}
    which implies that $a_3 = a_6$. Therefore, comparing \eqref{eq:al x>0} with \eqref{eq:al x<0} allows us to conclude that the constant coefficients in the domains $(-\infty, 0)$ and $(0, \infty)$ are equal for $f_X$.

    In section \ding{194}, where $x < y < 0$, the transformation simplifies to
    \[
        (u,v) = F_2(x,y) = (x-y,-y)
    \]
    which has Jacobian determinant $J = \begin{vmatrix} 1 & -1 \\ 0 & -1 \end{vmatrix} = -1$.
    Plugging this transformation into \eqref{eq:al1}, we get
    \begin{equation}\label{eq:al17}
        f_U(x-y) f_V(-y) = f_X(x) f_Y(y) \qquad \text{ for all } x<y<0.
    \end{equation}
    Substituting \eqref{eq:al u<0}, \eqref{eq:al v>0}, \eqref{eq:al x<0}, and \eqref{eq:al y<0} into \eqref{eq:al17} gives
    \begin{equation*}
        e^{a_4 (x-y) + a_5} e^{(a_4 + a_7)y + a_6 + a_9 - a_8} = e^{a_4 x + a_6} e^{a_7 y + a_9},
    \end{equation*}
    which implies that $a_5 = a_8$. Therefore, comparing \eqref{eq:al u>0} with \eqref{eq:al u<0} allows us to conclude that the constant coefficients in the domains $(-\infty, 0)$ and $(0, \infty)$ are equal for $f_U$.
    
    In section \ding{196}, where $x > 0$ and $y < 0$, the transformation simplifies to
    \[
        (u,v) = F_2(x,y) = (-y,-x)
    \]
    which has Jacobian determinant $J = \begin{vmatrix} 0 & -1 \\ -1 & 0 \end{vmatrix} = -1$.
    Plugging this transformation into \eqref{eq:al1}, we get
    \begin{equation}\label{eq:al18}
        f_U(-y) f_V(-x) = f_X(x) f_Y(y)\qquad \text{ for } x>0 \text{ and }y<0.
    \end{equation}
    Substituting \eqref{eq:al u>0}, \eqref{eq:al v<0}, \eqref{eq:al x>0}, and \eqref{eq:al y<0} into \eqref{eq:al18} gives
    \begin{equation}\label{eq:al19}
        e^{a_7 y + a_8} e^{-a_1 x + a_2} = e^{-a_1 x + a_3} e^{a_7 y + a_9}.
    \end{equation}
    Simplifying \eqref{eq:al19} yields
    \begin{equation*}
        a_8 + a_2 = a_3 + a_9.
    \end{equation*}
    Combining this with the fact that $a_3 = a_6$ and $a_5 = a_8$ implies that $a_9 = a_2 + a_5 - a_3$ and $a_2 = a_6 + a_9 - a_8$. Therefore, comparing \eqref{eq:al y>0} with \eqref{eq:al y<0} and \eqref{eq:al v>0} with \eqref{eq:al v<0} allows us to conclude that the constant coefficients in the domains $(-\infty, 0)$ and $(0, \infty)$ are equal for $f_Y$ and $f_V$.
    
    Finally, all that remains is to prove that the constant coefficients in their p.d.f.'s have the conjectured form. Using the fact that the p.d.f. of a random variable must integrate to 1, we obtain the densities of $X,Y,U,V$. Specifically, we combine $\eqref{eq:al x>0}$ and $\eqref{eq:al x<0}$ for $f_X$, $\eqref{eq:al y>0}$ and $\eqref{eq:al y<0}$ for $f_Y$, $\eqref{eq:al u>0}$ and $\eqref{eq:al u<0}$ for $f_U$, and $\eqref{eq:al v>0}$ and $\eqref{eq:al v<0}$ for $f_V$ to obtain
    \begin{align*}
        f_X(x) & = \frac{a_1a_4}{a_1 + a_4} (e^{-a_1x} \mathbf{1}_{[0,\infty)}(x) + e^{a_4x} \mathbf{1}_{(-\infty,0)}(x)), \quad & x \in \mathbb{R}, \\[5pt]
        f_Y(y) & = \frac{(a_1 + a_4) a_7}{a_1 + a_4 + a_7} (e^{-(a_1 + a_4)y} \mathbf{1}_{[0,\infty)}(y) + e^{a_7y} \mathbf{1}_{(-\infty,0)}(y)), \quad & y \in \mathbb{R}, \\[5pt]
        f_U(u) & = \frac{a_7a_4}{a_7 + a_4} (e^{-a_7u} \mathbf{1}_{[0,\infty)}(u) + e^{a_4u} \mathbf{1}_{(-\infty,0)}(u)), \quad & u \in \mathbb{R}, \\[5pt]
        f_V(v) & = \frac{(a_4 + a_7) a_1}{a_4 + a_7 + a_1} (e^{-(a_4 + a_7)v} \mathbf{1}_{[0,\infty)}(v) + e^{a_1v} \mathbf{1}_{(-\infty,0)}(v)), \quad & v \in \mathbb{R},
    \end{align*}
    where $a_1, a_4, a_7$ are positive constants (due to the fact that the densities are integrable). Replacing $a_1,a_4,a_7$ by $p,q,r$ gives us the p.d.f.'s that the conjecture predicts. We have thus proven the \textit{only if} part of the conjecture.
\end{proof}

\section{Proof of Theorem 1.3}
Similar to the proof of Theorem 1.2, we begin by partitioning $\mathbb{R}^2$ into sections so that the minimum functions in \eqref{eq:sExptrans} may be simplified. Define connected open subsets of $\mathbb{R}^2$ by
\begin{align*}
    \text{\ding{192}} :=\{{(x,y): x > -y}\}, \quad \text{\ding{193}} :=\{{(x,y): x < -y}\}.
\end{align*}
See Figure 2 below.

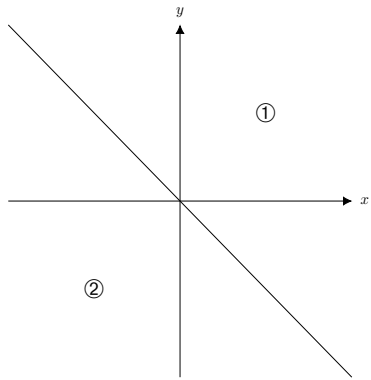
\begin{figure}[ht]
\centering
\resizebox{5cm}{5cm}{
    \begin{tikzpicture}
        \draw [{Latex[length = 2mm, width = 2mm]}-](0,4) -- (0,-4);
        \draw [{Latex[length = 2mm, width = 2mm]}-](4,0) -- (-4,0);
        \draw (-4,4) -- (4,-4);
        \node at (2,2) {\huge{\ding{192}}};
        \node at (-2,-2) {\huge{\ding{193}}};
        \node at (4.3,0) {$x$};
        \node at (0,4.3) {$y$};
    \end{tikzpicture}
}
\caption{Partition of $\mathbb{R}^2$}
\label{fig:2}
\end{figure}

Recall that the transformation in this theorem is
\[
    F_3(x,y) = (\min\{-x,y\}, y + x - \min\{-x,y\}) =: (u,v),
\]
which maps \ding{192} $\leftrightarrow$ \ding{192} and \ding{193} $\leftrightarrow$ \ding{193}. Therefore, defining the set $O := $ \ding{192} $\cup$ \ding{193}, we see that $F_3$ is a smooth involution on the open set $O$. Notice that $\mathbb{R}^2 \setminus O$ has Lebesgue measure $0$ and therefore has no effect on the distributions of the pairs of random variables $(X,Y)$ and $(U,V)$.  We will therefore only focus on the joint densities with inputs $(x,y)\in O$, which in turn implies $(u,v)=F_3(x,y)\in O$. 

Our approach again begins with the joint densities of $(U,V)$ and $(X,Y)$, which are related in the following manner
\begin{equation}\label{eq: invol density 2}
    f_{UV}(u,v) = |J| f_{XY}(x,y)   \qquad \text{ for } (x,y) \in O
\end{equation}
where $f_A(a)$ denotes the probability distribution function of $A$, and $J$ is the Jacobian determinant of the involution $F_3$. One can  check (as we do later) that $J\equiv -1$ throughout \ding{192}, \ding{193}, and therefore O.

\begin{proof}
    The \textit{if} part is easy to check by routine calculation using \eqref{eq: invol density 2}. We now prove the \textit{only if} part. Since $U,V$ are mutually independent and $X,Y$ are mutually independent, \eqref{eq: invol density 2} simplifies to
    \begin{equation}\label{eq:sExp1}
        f_U(u) f_V(v) = f_X(x) f_Y(y) \qquad \text{ for } (x,y) \in O.
    \end{equation}
    Since $f_X$ does not vanish on $[-c_1, c_2]$ and $f_Y$ does not vanish on $[-c_2, \infty)$ and $F_3: [-c_1, c_2] \times [-c_2, \infty) \rightarrow [-c_2, c_1] \times [-c_1, \infty)$ is a bijection and an involution on $\mathbb{R}^2$, \eqref{eq:sExp1} implies that $f_U$ does not vanish on $[-c_2, c_1]$ and $f_V$ does not vanish on $[-c_1, \infty)$. Since $F_3$ is a smooth involution on $O$, and $f_X$ and $f_Y$ are both twice-differentiable inside their respective supports, we know that $f_U$ and $f_V$ are twice differentiable on $[-c_2,c_1]$ and $[-c_1,\infty)$.
    
    In section \ding{192}, where $x > -y$, the transformation simplifies to
    \[
        (u,v) = F_3(x,y) = (-x,y+2x)
    \]
    whose Jacobian determinant is $J = \begin{vmatrix} -1 & 0 \\ 2 & 1 \end{vmatrix} = -1$.
    Plugging this transformation into \eqref{eq:sExp1} we get
    \begin{equation}\label{eq:sExp2}
        f_U(-x) f_V(y+2x) = f_X(x) f_Y(y) \qquad \text{ for all } x > -y.
    \end{equation}
    Restricting the domains of $x$ and $y$, taking logarithms of both sides of \eqref{eq:sExp2}, and setting $r_A(a) = \log f_A(a)$ gives
    \begin{equation}\label{eq:sExp3}
        r_U(-x) + r_V(y+2x) = r_X(x) + r_Y(y)
    \end{equation}
    for all $(x,y) \in [-c_1,c_2] \times [-c_2,\infty)$ such that $x > -y$.
    Differentiating \eqref{eq:sExp3} with respect to $y$ yields
    \begin{equation}\label{eq:sExp4}
        {r_V}'(y+2x) = {r_Y}'(y).
    \end{equation}
    Now differentiate \eqref{eq:sExp4} with respect to $x$ to obtain
    \begin{equation}\label{eq:sExp5}
        0 = 2{r_V}''(y+2x) = 2{r_V}''(v).
    \end{equation}
    The second equality holds for all $v>-c_1$ because in section \ding{192}, $y + 2x = v$. Thus $r_V$ must be linear on $(-c_1,\infty)$ and there must exist some real constants $a_1, a_2$ such that:
    \begin{equation}\label{eq:sExp6}
        r_V(v) = a_1 v + a_2, \qquad \text{ for all } v \in (-c_1, \infty).
    \end{equation}
    This implies that
    \begin{equation}\label{eq:sExp v}
        f_V(v) = e^{a_1 v + a_2} \qquad \text{ for all } v \in (-c_1, \infty).
    \end{equation}
    Substituting the derivative of \eqref{eq:sExp6} into \eqref{eq:sExp4} gives ${r_Y}'(y) = a_1$ for all $y \in (-c_2, \infty)$, which in turn implies the existence of some real constant $a_3$ such that 
    \begin{equation*}
        r_Y(y) = a_1 y + a_3 \qquad \text{ for all } y \in (-c_2, \infty),
    \end{equation*}
    which implies
    \begin{equation}\label{eq:sExp y}
        f_Y(y) = e^{a_1 y + a_3} \qquad \text{ for all } y \in (-c_2, \infty).
    \end{equation}
    We have thus shown that the p.d.f.'s of $Y$ and $V$ have exponential forms in $(-c_2, \infty)$ and $(-c_1, \infty)$.
    
    In section \ding{193}, where $x < -y$, the transformation simplifies to
    \[
        (u,v) = F_3(x,y) = (y,x)
    \]
    which has Jacobian determinant $J = \begin{vmatrix} 0 & 1 \\ 1 & 0 \end{vmatrix} = -1$.
    Plugging this into \eqref{eq:sExp1} we get
    \begin{equation}\label{eq:sExp7}
        f_U(y) f_V(x) = f_X(x) f_Y(y) \qquad \text{ for all } x < -y.
    \end{equation}
    Restricting the domains of $x$ and $y$, substituting \eqref{eq:sExp v} and \eqref{eq:sExp y} into \eqref{eq:sExp7}, taking the logarithm of both sides, and setting $r_A(a) = \log f_A(a)$ gives
    \begin{equation}\label{eq:sExp8}
        r_U(y) + a_1 x + a_2 = r_X(x) + a_1 y + a_3
    \end{equation}
    for all $(x,y) \in [-c_1,c_2] \times [-c_2,\infty)$ such that $x < -y$.
    Differentiating \eqref{eq:sExp8} with respect to $x$ and $y$ yields
    \begin{equation*}
        {r_X}'(x) = a_1 \qquad\text{ and } \qquad{r_U}'(y) = {r_U}'(u) = a_1.
    \end{equation*}
    The equality $r_U'(y) = r_U'(u)$ holds for all $u \in (-c_2, c_1)$ because in section \ding{193}, $y = u$. This implies that for some real constants $a_4, a_5$
    \begin{equation}\label{eq:sExp x}
        f_X(x) = e^{a_1 x + a_4} \qquad \text{ for all } x \in (-c_1, c_2),
    \end{equation}
    \begin{equation}\label{eq:sExp u}
        f_U(u) = e^{a_1 u + a_5} \qquad \text{ for all } u \in (-c_2, c_1).
    \end{equation}
    We have thus shown that the p.d.f.'s of $X$ and $U$ have exponential forms in $(-c_1, c_2)$ and $(-c_2, c_1)$.
    
    Recall that we showed $f_Y(y) = e^{a_1 y + a_3}$ on $(-c_2, \infty)$. Since a proper density function must integrate to 1, we must have $a_1 < 0$. Setting $\lambda = -a_1$, we see that $\lambda > 0$ and hence all four probability density functions have the desired rate parameter.
    
    We finish the proof by showing that $f_U$ and $f_V$ have supports corresponding to stExp and sExp distributions. Recall that we have already shown $f_U(u) > 0$ for all $u \in (-c_2, c_1)$ and $f_V(v) > 0$ for all $v \in (-c_1, \infty)$. All that remains is to show $f_U(u) = 0$ for all $u \in (-\infty, -c_2) \cup (c_1, \infty)$ and $f_V(v) = 0$ for all $v < -c_1$.
    
    By assumption, $f_X(x) > 0$ when $x \in [-c_1, c_2]$, $f_X(x) = 0$ otherwise, $f_Y(y) > 0$ when $y \in [-c_2, \infty)$ and $f_Y(y) = 0$ otherwise.
    
    
    
    By equation \eqref{eq:sExp2}, we have $f_U(-x) f_V(y+2x) = f_X(x) f_Y(y)$ throughout section \ding{192}, where $x > -y$. For all $x \in (-\infty, -c_1) \cup (c_2, \infty)$, $f_X(x) = 0$, so
    \begin{equation}\label{eq: end 1}
        f_U(-x) f_V(y+2x) = 0
    \end{equation}
    whenever $x \notin [-c_1,c_2]$ and $x > -y$. Taking any $y > \max\{-2x, 0\}$ gives $f_V(y+2x) > 0$ since $y+2x > 0 > -c_1$. Thus, \eqref{eq: end 1} implies $f_U(u) = 0$ for all $u \in (-\infty, -c_2) \cup (c_1, \infty)$ since $u = -x$.
    
    By equation \eqref{eq:sExp7}, we have $f_U(y) f_V(x) = f_X(x) f_Y(y)$ throughout section \ding{193}, where $x < -y$. For all $x < -c_1$, $f_X(x) = 0$, so
    \begin{equation}\label{eq: end 2}
        f_U(y) f_V(x) = 0
    \end{equation}
    whenever $x < -c_1 $ and $x < -y$. Taking any $y \in (-c_2, c_1)$ gives $f_U(y) > 0$. Thus, \eqref{eq: end 2} implies $f_V(v) = 0$ for all $v < -c_1$ since $v = x$.
    
    Combining all the above gives us the p.d.f.'s that the conjecture predicts. We have thus proven the \textit{only if} part of the conjecture.
\end{proof}

\section{Proof of Theorem 1.1}

Recall the transformation
\[
F_1(x,y) = \left(\dfrac{y(1 + \beta xy)}{1 + \alpha xy}, \dfrac{x(1 + \alpha xy)}{1 + \beta xy}\right) =: (u,v).
\]
Before we begin the proof, we make the useful observation that $uv = xy$. We also compute the following partial derivatives:
\begin{equation}
\begin{aligned}\label{partials}
    &\frac{\partial u}{\partial x} = \frac{(\beta - \alpha)y^2}{(1 + \alpha xy)^2} = \frac{(\beta - \alpha)u^2}{(1 + \beta uv)^2} & &\frac{\partial v}{\partial x} = \frac{1 + 2\alpha xy + \alpha \beta x^2y^2}{(1 + \beta xy)^2}  \\
    &\frac{\partial u}{\partial y} = \frac{1 + 2\beta xy + \alpha \beta x^2y^2}{(1 + \alpha xy)^2} & &\frac{\partial v}{\partial y} = \frac{(\alpha - \beta)x^2}{(1 + \beta xy)^2} = \frac{(\alpha - \beta)v^2}{(1 + \alpha uv)^2} \\
    &\frac{\partial^2 u}{\partial y \partial x} = \frac{2(\beta - \alpha)y}{(1 + \alpha xy)^3} = \frac{2(\beta - \alpha)u}{(1 + \alpha uv)^2(1 + \beta uv)} & &\frac{\partial^2 v}{\partial y \partial x} = \frac{2(\alpha - \beta)x}{(1 + \beta xy)^3} = \frac{2(\alpha - \beta)v}{(1 + \beta uv)^2(1 + \alpha uv)}.
    \end{aligned}
\end{equation}

Our approach again begins by analyzing the relationship between the joint densities of $(U,V)$ and $(X,Y)$, which are related in the following manner
\begin{equation}\label{eq:4.2}
    f_{UV}(u,v) = |J| f_{XY}(x,y) \qquad \text{ for all } x,y>0
\end{equation}
where $J$ is the Jacobian determinant of the transformation $F_1$.  Using \eqref{partials}  one can compute $J \equiv -1$.

\begin{proof}[Proof of Theorem 1.3]
    The \textit{if} part is easy to check by routine calculation using \eqref{eq:4.2}. We now prove the \textit{only if} part. Since $U,V$ are mutually independent and $X,Y$ are mutually independent, \eqref{eq:4.2} simplifies to
    \begin{equation}\label{eq:GIG trans}
        f_U(u) f_V(v) = f_X(x) f_Y(y).
    \end{equation}
    Since neither $f_X$ nor $f_Y$ vanish on $(0, \infty)$, \eqref{eq:GIG trans} implies that $f_U$ and $f_V$ do not vanish on $(0, \infty)$ as well. Since $F_1$ is a smooth involution, and $f_X$ and $f_Y$ are both twice differentiable, so are $f_U$ and $f_V$.
    Now taking logarithms and setting $r_A(a) = \log f_A(a)$ gives us
    \begin{equation}\label{eq:4.3}
        r_U(u) + r_V(v) = r_X(x) + r_Y(y).
    \end{equation}
    Now, taking mixed partials $\frac{\partial^2}{\partial y \partial x}$ of both sides we get
    \begin{align*}
        0 &= \frac{\partial^2}{\partial y \partial x} r_U(u) + \frac{\partial^2}{\partial y \partial x} r_V(v) \\
        &= \frac{\partial}{\partial y} \left( \frac{\partial u}{\partial x} r_U{'}(u) \right) + \frac{\partial}{\partial y} \left( \frac{\partial v}{\partial x} r_V{'}(v) \right) \\
        &= \frac{\partial u}{\partial y} \frac{\partial u}{\partial x} r_U{''}(u) + \frac{\partial^2 u}{\partial y \partial x} r_U{'}(u) + \frac{\partial v}{\partial y} \frac{\partial v}{\partial x} r_V{''}(v) + \frac{\partial^2 v}{\partial y \partial x} r_V{'}(v).
    \end{align*}
    Substituting in the values of the partials from \eqref{partials}, using the fact that $xy = uv$, multiplying both sides by the common denominator $\frac{(1 + \beta uv)^2(1 + \alpha uv)^2}{\beta - \alpha}$, and rearranging gives us
    \begin{equation}\label{eq:4.4}
        \begin{aligned}
            &(1 + 2\beta uv + \alpha \beta u^2v^2) u^2 r_U{''}(u) + (1 + \beta uv) 2u r_U{'}(u) \\ =& (1 + 2\alpha uv + \alpha \beta u^2v^2) v^2 r_V{''}(v) + (1 + \alpha uv) 2v r_V{'}(v).
        \end{aligned}
    \end{equation}
    
    Since $(u,v) = F_1(x,y)$ is an involution, equation \eqref{eq:4.4} holds for all $(u,v) \in (0,\infty)^2$. By first fixing any $v>0$, we see the right-hand side (and therefore the left-hand side) of \eqref{eq:4.4} must be a second degree polynomial in $u$ for all $u>0$.  Since second degree polynomials are smooth, it follows that $r_U$ is smooth on $(0,\infty)$.   Similarly, for each fixed $u>0$, the left-hand side (and therefore the right-hand side) must be a second degree polynomial in $v$ for all $v>0$.  It now follows that $r_V$ is smooth on $(0,\infty)$.  The fact that $r_V$ and $r_U$ are both smooth now implies that when writing the right-hand side of \eqref{eq:4.4} out as a second degree polynomial in $u$, each of the coefficients has a smooth dependence on $v$. This implies that each of these coefficient must in turn themselves be at most second degree polynomials in $v$ with no $u$-dependence.  We may therefore set both sides of \eqref{eq:4.4} equal to
     \[
         a_0 + a_1 u + a_2 v + a_3 u^2 + a_4 uv + a_5 v^2 + a_6 u^2 v + a_7 u v^2 + a_8 u^2v^2,
     \]
     for $u,v>0$ where $a_0, a_1, ... ,a_8 \in \mathbb{R}$ are fixed constants.
    
    Taking the limits of equation \eqref{eq:4.4} and the polynomial as $v \rightarrow 0^+$ and $u \rightarrow 0^+$, we get respectively
    \begin{align*}
             (u^2 r_U{'}(u)){'} &=u^2 r_U{''}(u) + 2u r_U{'}(u) = a_0 + a_1 u + a_3 u^2 \\
            (v^2 r_V{'}(v)){'} &=v^2 r_V{''}(v) + 2v r_V{'}(v) = a_0 + a_2 v + a_5 v^2.
    \end{align*}
    These two differential equations hold for all $u,v>0$ and can be explicitly solved for resulting in
    \begin{equation}
    \begin{aligned}\label{eq:4.5}
        r_U(u) &= a_0 \ln u + \frac{a_1}{2} u + \frac{a_3}{6} u^2 - \frac{b_0}{u} + b_2 \\
        r_V(v) &= a_0 \ln v + \frac{a_2}{2} v + \frac{a_5}{6} v^2 - \frac{b_1}{v} + b_3
    \end{aligned}
    \end{equation}
    where $b_0, b_1, b_2, b_3 \in \mathbb{R}$ are fixed constants.
    Notice that the forms of $r_U$ and $r_V$ are very close to that of the logarithm of the p.d.f.\ of the GIG distribution. To finish the proof, we just have to show that $a_3 = a_5 = 0$ and that the remaining parameters are related exactly as conjectured.
    
    Substituting the derivatives of $r_U$ and $r_V$ into equation \eqref{eq:4.4} gives us
    \begin{align*}
        &(1 + 2\beta uv + \alpha \beta u^2v^2) \left( -a_0 + \dfrac{a_3}{3} u^2 - \dfrac{2b_0}{u} \right) + 2(1 + \beta uv) \left( a_0 + \dfrac{a_1}{2} u + \dfrac{a_3}{3} u^2 + \dfrac{b_0}{u} \right) \\ =& (1 + 2\alpha uv + \alpha \beta u^2v^2) \left( -a_0 + \dfrac{a_5}{3} v^2 - \dfrac{2b_1}{v} \right) + 2(1 + \alpha uv) \left( a_0 + \dfrac{a_2}{2} v + \dfrac{a_5}{3} v^2 + \dfrac{b_1}{v} \right).
    \end{align*}
    Combining like terms and simplifying, we have
    \begin{equation}\label{eq:4.6}
        \begin{aligned}
            &(3 + 4\beta uv + \alpha \beta u^2v^2) \dfrac{a_3}{3} u^2 + (1 + \beta uv) a_1 u - 2b_0 (\beta v + \alpha \beta uv^2) \\ =& (3 + 4\alpha uv + \alpha \beta u^2v^2) \dfrac{a_5}{3} v^2 + (1 + \alpha uv) a_2 v - 2b_1 (\alpha u + \alpha \beta u^2v).
        \end{aligned}
    \end{equation}
   Both sides of \eqref{eq:4.6} are polynomials in $u$ and $v$, so we equate like terms.  Since the left-hand side has no $v^2$ term and the right-hand side has no $u^2$ term, it follows that $a_3=a_5=0$. Continuing to equate the $u$ and $v$ terms we get
   \[
   a_1=-2b_1 \alpha\qquad \text{ and } \qquad a_2=-2 b_0 \beta.
   \]
   Now set $c_1=b_0$, $c_2=b_1$, $\lambda=a_0+1$, $Z_U=e^{-b_2}$, $Z_V=e^{-b_3}$ and plug these into \eqref{eq:4.5} to get

    \begin{equation*}
        \begin{dcases}
            r_U(u) = (\lambda-1) \ln u - (c_2 \alpha u + \frac{c_1}{u}) + b_2 \\
            r_V(v) = (\lambda-1) \ln v - (c_1 \beta v + \frac{c_2}{v}) + b_3
        \end{dcases}
        \implies
       \begin{dcases}
            f_U(u) = \dfrac{1}{Z_U} u^{\lambda - 1} e^{-\left( c_2 \alpha u + \frac{c_1}{u} \right)} \\
            f_V(v) = \dfrac{1}{Z_V} v^{\lambda - 1} e^{-\left( c_1 \beta v + \frac{c_2}{v} \right)}.
        \end{dcases}
    \end{equation*}
    The fact that $f_U$ and $f_V$ are densities, and are therefore integrable, implies the constants $c_1, c_2$ have the right signs and that 
     \[
         U \sim \textup{GIG}(\lambda,c_2\alpha,c_1) \qquad \text{ and } \qquad  V \sim \textup{GIG}(\lambda,c_1\beta,c_2).
    \]
    Finally, since $(U,V)=F_1(X,Y)$ and $F_1$ is an involution, we have 
    $(X,Y)=F_1(U,V)$.  Now the \textit{if} part of the theorem implies that $X$ and $Y$ have GIG distributions with the desired parameters.

\end{proof}

\bibliographystyle{amsplain}

\end{document}